# NECESSARY AND SUFFICIENT CONDITION FOR A SET OF MATRICES TO COMMUTE

M. De la Sen

Department of Electricity and Electronics. **Faculty of Science and Technology**.

Campus of Leioa (Bizkaia). Aptdo. 644- Bilbao, SPAIN

**Abstract**. This paper investigates the necessary and sufficient condition for a set of (real or complex) matrices to commute. It is proved that the commutator $[A,B]=0$ for two matrices A and B if and only if a vector v (B) defined uniquely from the matrix B is in the null space of a well- structured matrix defined as the Konecker sum $A \oplus (-A^*)$, which is always rank defective. This result is extendable directly to any countable set of commuting matrices. Complementary results are derived concerning commutators of certain matrices with functions of matrices $f(A)$ which extends the well- known sufficiency – type commuting result $[A, f(A)] = 0$.

**Keywords**: Commuting matrices, Jordan canonical form, Kronecker product, Kronecker sum, Function of a matrix.

### 1. Introduction

The problem of commuting matrices is very relevant in certain problems of Engineering and Physics . In particular, such a problem is of crucial interest related to discrete Fourier transforms, normal modes in dynamic systems or commuting matrices dependent on a parameter (see, for instance, [1-3]). It is well-known that commuting matrices have at least a common eigenvector and also, a common generalized eigenspace, [4-5]. A less restrictive problem of interest in the above context is that of almost commuting matrices, roughly speaking, the norm of the commutator is sufficiently small, [5-6]. A very relevant related result is that the sum of matrices which commute is an infinitesimal generator of a $C_0$-semigroup. This leads to a well-known result in Systems Theory establishing that that the matrix function $e^{A_1 t_1 + A_2 t_2} = e^{A_1 t_1} e^{A_2 t_2}$ is a fundamental (or state transition) matrix for the cascade of the time invariant differential systems $\dot{x}_1(t) = A_1 x_1(t)$, operating on a time $t_1$, and $\dot{x}_2(t) = A_2 x_2(t)$, operating on a time $t_2$, provided that $A_1$ and $A_2$ commute (see, for instance, [7-11] ).The problem of commuting matrices is also of relevant interest in dynamic switched systems, namely, those which possess several parameterizations one of each is activated at each current time interval. If the matrices of dynamics of all the parameterizations commute then there exists a common Lyapunov function for all those parameterizations and any arbitrary switching rule operating at any time instant maintains the global stability of the switched rule provided that all the parameterizations are stable, [7]. However, in the case that there is no common Lyapunov function for all the parameterizations , a minimum residence time at each active parameterization is needed to maintain the global stability of the switched system so that the switching rule among distinct parameterizations is not arbitrary, [12-13]. This fact implies that . Parallel results apply for switched time-delay systems subject to point delays under zero or sufficiently small delays when the matrices defining the delay-free dynamics of the various parameterizations commute, [10-11].

This paper formulates the necessary and sufficient condition for any countable set of (real or complex) matrices to commute. The sequence of obtained the results is as follows. Firstly, the commutation of two real matrices is investigated in Section 2. The necessary and sufficient condition for two matrices to commute is that a vector defined uniquely from the entries of any of the two given matrices belongs to the null space of the Kronecker sum of the other matrix and its minus transpose. The above result allows a simple algebraic characterization and computation of the set of commuting matrices with a given one. It also exhibits counterparts for the necessary and sufficient condition for two matrices not to commute. The results are then extended to the necessary and sufficient condition for commutation of any set of real matrices in Section 3. In Section 4, the above results are directly extended to the case of complex matrices in two very simple ways, namely, either decomposition of the associated algebraic system of complex matrices into two real ones or by manipulating it directly as a complex algebraic system of equations. Basically, the results for the real case are directly extendable by replacing transposes by conjugate transposes. Finally, further results concerning the commutators of matrices with matrix functions are also discussed in Section 4.

### 1.1. Notation

$[A, B]$ is the commutator of the square matrices A and B.

$A \otimes B := (a_{ij} B)$ is the Kronecker (or direct) product of $A := (a_{ij})$ and B.

$A \oplus B := A \otimes I_n + I_n \otimes B$ is the Kronecker sum of the square matrices $A := (a_{ij})$ and both of order n, where $I_n$ is the n-th identity matrix.

$A^T$ is the transpose of the matrix A and $A^*$ is the conjugate transpose of the complex matrix A. For any matrix A, Im A and Ker A are its associate range (or image) subspace and null space, respectively. Also, rank (A) is the rank of A which is the dimension of Im (A) and det (A) is the determinant of the square matrix A.

$v(A) = (a_1^T, a_2^T, ...., a_n^T)^T \in \mathbb{C}^{n^2}$ if $a_i^T := (a_{i1}, a_{i2}, ...., a_{in})$ is the i-th row of the square matrix A.

$\sigma(A)$ is the spectrum of A; $\bar{n} := \{1, 2, ..., n\}$. If $\lambda_i \in \sigma(A)$ then there exist positive integers $\mu_i$ and $\nu_i \leq \mu_i$ which are, respectively, its algebraic and geometric multiplicity; i.e. the times it is repeated in the characteristic polynomial of A and the number of its associate Jordan blocks, respectively. The integer $\mu \leq n$ is the number of distinct eigenvalues and the integer $m_i$, subject to $1 \leq m_i \leq \mu_i$, is the index of $\lambda_i \in \sigma(A)$; $\forall i \in \bar{\mu}$, that is, its multiplicity in the minimal polynomial of A.

$A \sim B$ denotes a similarity transformation from A to $B = T^{-1} A T$ for given $A, B \in \mathbb{R}^{n \times n}$ for some non-singular $T \in \mathbb{R}^{n \times n}$. $A \approx B = E A F$ means that there is an equivalence transformation for given $A, B \in \mathbb{R}^{n \times n}$ for some non-singular $E, F \in \mathbb{R}^{n \times n}$.

A linear transformation from $\mathbb{R}^n$ to $\mathbb{R}^n$, represented by the matrix $T \in \mathbb{R}^{n \times n}$, is denoted identically to such a matrix in order to simplify the notation. If $V \neq \text{Dom} T \equiv \mathbb{R}^n$ is a subspace of $\mathbb{R}^n$ then $\operatorname{Im} T(V) := \{ Tz : z \in V \}$ and $\operatorname{Ker} T(V) := \{ z \in V : Tz = 0 \in \mathbf{R}^n \}$. If $V \equiv \mathbf{R}^n$, the notation is simplified to $\operatorname{Im} T := \{ Tz : z \in \mathbf{R}^n \}$ and $\operatorname{Ker} T := \{ z \in \mathbf{R}^n : Tz = 0 \in \mathbf{R}^n \}$.

The symbols "$\wedge$" and "$\vee$" stand for logic conjunction and disjunction, respectively.

The abbreviation "iff" stands for "if and only if".

The notation card U stands for the cardinal of the set U.

## 2. Results concerning the sets of commuting and no commuting matrices with a given one

Consider the sets $C_A := \{ X \in \mathbf{R}^{n \times n} : [A, X] = 0 \} \neq \varnothing$, of matrices which commute with A, and $\overline{C}_A := \{ X \in \mathbf{R}^{n \times n} : [A, X] \neq 0 \}$, of matrices which do not commute with A; $\forall A \in \mathbf{R}^{n \times n}$. Note that $0 \in \mathbf{R}^{n \times n} \cap C_A$; i.e. the zero n- matrix commutes with any n-matrix so that, equivalently, $0 \notin \mathbf{R}^{n \times n} \cap \overline{C}_A$ and then $C_A \cap \overline{C}_A = \varnothing$; $\forall A \in \mathbf{R}^{n \times n}$. The following two basic results follow concerning commutation and non- commutation of two matrices:

**Propositions 2.1**. (i) $C_A := \{ X \in \mathbf{R}^{n \times n} : v(X) \in \operatorname{Ker}(A \oplus (-A^T)) \}$, and equivalently,

$[A, X] = 0 \Leftrightarrow v(X) \in \operatorname{Ker}(A \oplus (-A^T))$.

(ii) $\overline{C}_A := \mathbf{R}^{n \times n} \setminus C_A = \{ X \in \mathbf{R}^{n \times n} : v(X) \notin \operatorname{Ker}(A \oplus (-A^T)) \} \equiv \{ X \in \mathbf{R}^{n \times n} : v(X) \in \operatorname{Im}(A \oplus (-A^T)) \}$, and, equivalently,

$[A, X] \neq 0 \Leftrightarrow v(X) \in \operatorname{Im}(A \oplus (-A^T))$.

(iii) $B \in C_A := \{ X \in \mathbf{R}^{n \times n} : v(X) \in \operatorname{Ker}(A \oplus (-A^T)) \}$

$\Leftrightarrow A \in C_B := \{ X \in \mathbf{R}^{n \times n} : v(X) \in \operatorname{Ker}(B \oplus (-B^T)) \}$

**Proof:** (i) –(ii) First note by inspection that $\varnothing \neq C_A \supset \{0, A\}$; $\forall A \in \mathbf{R}^{n \times n}$. Also,

$[A, X] = AX - XA = (A \otimes I_n - I_n \otimes A^T) v(X) = (A \oplus (-A^T)) v(X) = 0$

$\Rightarrow v(X) \in \operatorname{Ker}(A \oplus (-A^T))$ and Propositions 2.1(i)-(ii) have been proved since there is an isomorphism $f : \mathbf{R}^{n^2} \leftrightarrow \mathbf{R}^{n \times n}$ defined by $f(v(X)) = X$; $\forall X \in \mathbf{R}^{n \times n}$ for $v(X) = (x_1^T, x_2^T, \ldots, x_n^T)^T \in \mathbf{R}^{n^2}$ if $x_i^T := (x_{i1}, x_{i2}, \ldots, x_{in})$ is the i-th row of the square matrix X.

(iii) It is a direct consequence of Proposition 2.1 (iii) and the symmetry property of the commutator of two commuting matrices $B \in C_A \Leftrightarrow [A, B] = [B, A] = 0 \Leftrightarrow A \in C_B$. □

**Proposition 2.2**.

$\operatorname{rank}(A \oplus (-A^T)) < n^2 \Leftrightarrow \operatorname{Ker}(A \oplus (-A^T)) \neq 0 \Leftrightarrow 0 \in \sigma(A \oplus (-A^T)) \Leftrightarrow \exists X (\neq 0) \in C_A$; $\forall A \in \mathbf{R}^{n \times n}$.

**Proof**: $[A, A] = 0$; $\forall A \in \mathbf{R}^{n \times n} \Rightarrow \exists \mathbf{R}^{n^2} \ni 0 \neq v(A) \in \operatorname{Ker}(A \oplus (-A^T))$; $\forall A \in \mathbf{R}^{n \times n}$. As a result,

$\operatorname{Ker}(A \oplus (-A^T)) \neq 0 \in \mathbf{R}^{n^2}$; $\forall A \in \mathbf{R}^{n \times n} \Leftrightarrow \operatorname{rank}(A \oplus (-A^T)) < n^2$; $\forall A \in \mathbf{R}^{n \times n}$

so that $0 \in \sigma(A \oplus (-A^T))$.

Also, $\exists X(\neq 0) \in \mathbf{R}^{n \times n} : [A, X] = 0 \Leftrightarrow X \in C_A$ since $\mathrm{Ker}\left(A \oplus \left(-A^T\right)\right) \neq 0 \in \mathbf{R}^{n^2}$.

Then, Proposition 2.2 has been proved. □

The subsequent result is stronger than Proposition 2.2.

**Theorem 2.3**. The following properties hold:

**(i)** The spectrum of $A \oplus \left(-A^T\right)$ is $\sigma\left(A \oplus \left(-A^T\right)\right) = \{\bar{\lambda}_{ij} = \lambda_i - \lambda_j : \lambda_i, \lambda_j \in \sigma(A); \forall i, j \in \bar{n}\}$ and possesses $\bar{\nu}$ Jordan blocks in its Jordan canonical form of, subject to the constraints $n^2 \geq \bar{\nu} = \dim S = \left(\sum_{i=1}^{\mu} \nu_i\right)^2 \geq \bar{\nu}(0)$, and $0 \in \sigma\left(A \oplus \left(-A^T\right)\right)$ with an algebraic multiplicity $\bar{\mu}(0)$ and with a geometric multiplicity $\bar{\nu}(0)$ subject to the constraints:

$$n^2 = \left(\sum_{i=1}^{\mu} \mu_i\right)^2 \geq \bar{\mu}(0) \geq \sum_{i=1}^{\mu} \mu_i^2 \geq \bar{\nu}(0) = \sum_{i=1}^{\mu} \nu_i^2 \geq n \tag{2.1}$$

where:

a) $S := \mathrm{span}\{z_i \otimes x_j, \forall i, j \in \bar{n}\}$, $\mu_i = \mu(\lambda_i)$ and $\nu_i = \nu(\lambda_i)$ are, respectively, the algebraic and the geometric multiplicities of $\lambda_i \in \sigma(A)$, $\forall i \in \bar{n}$; $\mu \leq n$ is the number of distinct $\lambda_i \in \sigma(A)$ $(i \in \mu)$, $\bar{\mu}_i = \bar{\mu}(\bar{\lambda}_{ij})$ and $\bar{\nu}_{ij} = \bar{\nu}(\bar{\lambda}_{ij})$, are, respectively, the algebraic and the geometric multiplicity of $\bar{\lambda}_{ij} = (\lambda_i - \lambda_j) \in \sigma\left(A \oplus \left(-A^T\right)\right)$, $\forall i, j \in \bar{n}$; $\mu \leq n$, and

b) $x_j$ and $z_i$ are, respectively, the right eigenvectors of $A$ and $A^T$ with respective associated eigenvalues $\lambda_j$ and $\lambda_i$; $\forall i, j \in \bar{n}$.

**(ii)** $\dim \mathrm{Im}\left(A \oplus \left(-A^T\right)\right) = \mathrm{rank}\left(A \oplus \left(-A^T\right)\right) = n^2 - \bar{\nu}(0)$

$$\Leftrightarrow \dim \mathrm{Ker}\left(A \oplus \left(-A^T\right)\right) = \bar{\nu}(0); \forall A \in \mathbf{R}^{n \times n} \tag{2.2}$$

**Proof**: **(i)** Note that

$$\sigma(A) = \sigma\left(A^T\right) \Rightarrow \sigma\left(A \oplus \left(-A^T\right)\right) := \{\mathbf{C} \ni \eta = \lambda_k - \lambda_\ell; \forall \lambda_k, \lambda_\ell \in \sigma(A); \forall k, \ell \in \bar{n}\}$$

$$= \sigma_0\left(A \oplus \left(-A^T\right)\right) \cup \overline{\sigma_0\left(A \oplus \left(-A^T\right)\right)}$$

where

$$\sigma_0\left(A \oplus \left(-A^T\right)\right) = \{\lambda \in \sigma\left(A \oplus \left(-A^T\right)\right) : \lambda = 0\}$$

$$\overline{\sigma_0\left(A \oplus \left(-A^T\right)\right)} = \{\lambda \in \sigma\left(A \oplus \left(-A^T\right)\right) : \lambda \neq 0\} = \sigma\left(A \oplus \left(-A^T\right)\right) \setminus \sigma_0\left(A \oplus \left(-A^T\right)\right)$$

Furthermore, $\sigma\left(A \oplus \left(-A^T\right)\right) := \{\mathbf{C} \ni \lambda = \lambda_j - \lambda_i : \lambda_i, \lambda_j \in \sigma(A); \forall i, j \in \bar{n}\}$ and $z_i \otimes x_j$ is a right eigenvector of $A \oplus \left(-A^T\right)$ associated with its eigenvalue $\bar{\lambda}_{ji} = \lambda_j - \lambda_i$. $\lambda = \lambda_j - \lambda_i \in \sigma\left(A \oplus \left(-A^T\right)\right)$ has algebraic and geometric multiplicities $\bar{\mu}_{ji}$ and $\bar{\nu}_{ji}$, respectively;

$\forall i, j \in \bar{n}$, since $x_j$ and $z_i$ are, respectively, the right eigenvectors of $A$ and $A^T$ with associated eigenvalues $\lambda_j$ and $\lambda_i$; $\forall i, j \in \bar{n}$.

Let $J_A$ be the Jordan canonical form of $A$. It is first proved that there exists a non-singular $T \in \mathbf{R}^{n^2 \times n^2}$ such that $J_A \oplus (-J_{A^T}) = T^{-1}(A \oplus (-A^T))T$. The proof is made by direct verification by using the properties of the Kronecker product, with $T = P \otimes P^T$ for a non-singular $P \in \mathbf{R}^{n \times n}$ such that $A \sim J_A = P^{-1}AP$, as follows:

$$T^{-1}(A \oplus (-A^T))T = (P \otimes P^T)^{-1}(A \otimes I_n)(P \otimes P^T) - (P \otimes P^T)^{-1}(I_n \otimes A^T)(P \otimes P^T)$$
$$= (P^{-1}AP) \otimes (P^{-T}I_n P^T) - (P^{-1}I_n P) \otimes (P^{-T}A^T P^T)$$
$$= (P^{-1}AP) \otimes I_n - I_n \otimes (P^{-T}A^T P^T)$$
$$= J_A \otimes I_n - I_n \otimes J_{A^T} = J_A \otimes I_n + I_n \otimes (-J_{A^T}) = J_A \oplus (-J_{A^T})$$

and the result has been proved. Thus, $\text{rank}(A \oplus (-A^T)) = \text{rank}(J_A \oplus (-J_{A^T}))$. It turns out that $P$ is, furthermore, unique except for multiplication by any nonzero real constant. Otherwise, if $T \neq P \otimes P^T$, then there would exist a non-singular $Q \in \mathbf{R}^{n \times n}$ with $Q \neq \alpha I_n$; $\forall \alpha \in \mathbf{R}$ such that $T = Q(P \otimes P^T)^{-1}Q$ so that $T^{-1}(A \oplus (-A^T))T \neq J_A \oplus (-J_{A^T})$ provided that

$$(P \otimes P^T)^{-1}(A \oplus (-A^T))(P \otimes P^T) = J_A \oplus (-J_{A^T})$$

Thus, note that:

$$\text{card } \sigma(A \oplus (-A^T)) = n^2 = \sum_{i=1}^{\mu} \bar{\mu}_{ii} = \left(\sum_{i=1}^{\mu} \mu_i\right)^2 \geq \bar{\mu}(0) = \sum_{i=1}^{\mu} \bar{\mu}_{ii} = \sum_{i=1}^{\mu} \mu_i^2 \geq \bar{\nu}$$

$$\geq \bar{\nu}(0) = \sum_{i=1}^{\mu} \bar{\nu}_{ii} = \sum_{i=1}^{\mu} \nu_i^2 = \left(\sum_{i=1}^{\mu}\sum_{j=1}^{\mu} \bar{\nu}_{ij}\right)^2 - 2\sum_{i=1}^{\mu}\sum_{j(\neq i)=1}^{\mu} \bar{\nu}_{ij} = \bar{\nu} - 2\sum_{i=1}^{\mu}\sum_{j(\neq i)=1}^{\mu} \bar{\nu}_{ij} \geq n$$

Those results follow directly from the properties of the Kronecker sum $A \oplus B$ of n-square real matrices $A$ and $B = -A^T$ since direct inspection leads to:

(1) $0 \in \sigma(A \oplus (-A^T))$ with algebraic multiplicity $\bar{\mu}(0) \geq \sum_{i=1}^{\mu} \bar{\mu}_{ii} = \sum_{i=1}^{\mu} \mu_i^2 \geq \sum_{i=1}^{\mu} \nu_i^2 \geq n$ since there are at least $\sum_{i=1}^{n} \mu_i^2$ zeros in $\sigma(A \oplus (-A^T))$ (i.e. the algebraic multiplicity of $0 \in \sigma(A \oplus (-A^T))$ is at least $\sum_{i=1}^{n} \mu_i^2$) and since $\nu_i \geq 1$; $\forall i \in \bar{n}$. Also, a simple computation of the number of eigenvalues of $A \oplus (-A^T)$ yields $\text{card } \sigma(A \oplus (-A^T)) = n^2 = \sum_{i=1}^{\mu} \bar{\mu}_{ii} = \left(\sum_{i=1}^{\mu} \mu_i\right)^2$.

(2) The number of linearly independent vectors in S is $\bar{v} = \sum_{i=1}^{\mu}\sum_{j=1}^{\mu} \bar{v}_{ij} = \left(\sum_{i=1}^{\mu} v_i\right)^2 \geq \sum_{i=1}^{\mu} \bar{v}_{ii} = \sum_{i=1}^{\mu} v_i^2$ since the total number of Jordan blocks in the Jordan canonical form of A is $\sum_{i=1}^{\mu} v_i$.

(3) The number of Jordan blocks associated with $0 \in \sigma\left(A \oplus \left(-A^T\right)\right)$ in the Jordan canonical form of $\left(A \oplus \left(-A^T\right)\right)$ is $\bar{v}(0) = \sum_{i=1}^{\mu} v_i^2 \leq \bar{v}$, with $\bar{v}_{ii} = v_{ii}^2$; $\forall i \in \bar{n}$. Thus:

$$\text{card } \sigma_0\left(A \oplus \left(-A^T\right)\right) = \sum_{i=1}^{\mu} \bar{\mu}_{ii} = \sum_{i=1}^{\mu} \mu_i^2, \qquad \text{card } \bar{\sigma}_0\left(A \oplus \left(-A^T\right)\right) = n^2 - \sum_{i=1}^{\mu} \mu_i^2$$

$$\text{rank}\left(A \oplus \left(-A^T\right)\right) = n^2 - \bar{v}(0) = n^2 - \sum_{i=1}^{\mu} v_i^2, \quad \text{dim Ker}\left(A \oplus \left(-A^T\right)\right) = \bar{v}(0) = \sum_{i=1}^{\mu} v_i^2$$

(4) There are at least $\bar{v}(0)$ linearly independent vectors in $S := \text{span}\{z_i \otimes x_j, \forall i, j \in \bar{n}\}$. Also, the total number of Jordan blocks in the Jordan canonical form of $\left(A \oplus \left(-A^T\right)\right)$ is $\bar{v} = \dim S = \left(\sum_{i=1}^{\mu}\sum_{j=1}^{\mu} \bar{v}_{ij}\right) = \left(\sum_{i=1}^{\mu} v_i\right)^2 = \bar{v}(0) + 2\sum_{i=1}^{\mu}\sum_{j(\neq i)=1}^{\mu} \bar{v}_{ij} \geq \bar{v}(0)$.

Property (i) has been proved. Property **(ii)** follows directly from the orthogonality in $\mathbf{R}^{n^2}$ of its range and null subspaces. □

Expressions which calculate the sets of matrices which commute and which do not commute with a given one are obtained in the subsequent result:

**Theorem 2.4**. The following properties hold:

(i) $X \in C_A$ iff $\left(A \oplus \left(-A^T\right)\right) v(X) = 0 \Leftrightarrow X \in C_A$ iff $v(X) = -F\left(v^T(\bar{X}_2) \bar{A}_{12}^T \bar{A}_{11}^{-T}, v^T(\bar{X}_2)\right)^T$

for any $v(\bar{X}_2) \in \text{Ker}\left(\bar{A}_{22} - \bar{A}_{21} \bar{A}_{11}^{-1} \bar{A}_{12}\right)$, where $E, F \in \mathbf{R}^{n^2 \times n^2}$ are permutation matrices and $\bar{X} \in \mathbf{R}^{n \times n}$ and $v(\bar{X}) \in \mathbf{R}^{n^2}$ are defined such that:

(a) $v(\bar{X}) := F^{-1} v(X)$, $A \oplus \left(-A^T\right) \approx \bar{A} := E\left(A \oplus \left(-A^T\right)\right) F$; $\forall X \in C_A$ (2.3)

where $v(\bar{X}) = (v^T(\bar{X}_1), v^T(\bar{X}_2))^T \in \mathbf{R}^{n^2}$ with $v(\bar{X}_1) \in \mathbf{R}^{\bar{v}(0)}$ and $v(\bar{X}_2) \in \mathbf{R}^{n^2 - \bar{v}(0)}$

(b) The matrix $\bar{A}_{11} \in \mathbf{R}^{\bar{v}(0) \times \bar{v}(0)}$ is non-singular in the block matrix partition $\bar{A} := \text{Block matrix}\left(\bar{A}_{ij}; i, j \in \bar{2}\right)$ with $\bar{A}_{12} \in \mathbf{R}^{\bar{v}(0) \times n^2}$, $\bar{A}_{21} \in \mathbf{R}^{(n^2 - \bar{v}(0)) \times \bar{v}(0)}$ and $\bar{A}_{22} \in \mathbf{R}^{(n^2 - \bar{v}(0)) \times (n^2 - \bar{v}(0))}$.

**(ii)** $X \in \overline{C}_A$, for any given $A(\neq 0) \in \mathbf{R}^{n \times n}$, iff

$$\left(A \oplus \left(-A^T\right)\right) v(X) = v(M) \qquad (2.4)$$

for some $M(\neq 0) \in \mathbf{R}^{n \times n}$ such that:

$$\text{rank}\left(A \oplus \left(-A^T\right)\right) = \text{rank}\left(A \oplus \left(-A^T\right), v(M)\right) = n^2 - \overline{v}(0) \qquad (2.5)$$

Also,

$$\overline{C}_A := \{X \in \mathbf{R}^{n \times n} : \left(A \oplus \left(-A^T\right)\right) v(X) = v(M) \text{ for any } M(\neq 0) \in \mathbf{R}^{n \times n} \text{ satisfying}$$
$$\text{rank}\left(A \oplus \left(-A^T\right)\right) = \text{rank}\left(A \oplus \left(-A^T\right), v(M)\right) = n^2 - \overline{v}(0)\} \qquad (2.6)$$

Also, with the same definitions of $E$, $F$ and $\overline{X}$ in (i), $X \in \overline{C}_A$ iff

$$v(X) = F\left(v^T\left(\overline{M}_1\right)\overline{A}_{11}^{-T} - v^T\left(\overline{X}_2\right)\overline{A}_{12}^T \overline{A}_{11}^{-T}, v^T\left(\overline{X}_2\right)\right)^T \qquad (2.7)$$

where $v(\overline{X}_2)$ is any solution of the compatible algebraic system

$$\left(\overline{A}_{22} - \overline{A}_{21} \overline{A}_{11}^{-1} \overline{A}_{12}\right) v(\overline{X}_2) = v(\overline{M}_2) - \overline{A}_{21} \overline{A}_{11}^{-1} v(\overline{M}_1) \qquad (2.8)$$

for some $M(\neq 0) \in \mathbf{R}^{n \times n}$ such that $\overline{X}, \overline{M} \in \mathbf{R}^{n \times n}$ and are defined according to $v(X) = Fv(\overline{X})$ and $\overline{M} = EMF \approx M(\neq 0) \in \mathbf{R}^{n \times n}$ and $v(\overline{M}) = Ev(M) = E\left(v_1^T(M), v_2^T(M)\right)^T$.

**Proof**: First note from Proposition 2.1 that $X \in C_A$ iff $\left(A \oplus \left(-A^T\right)\right) v(X) = 0$ since $v(X) \in \text{Ker}\left(A \oplus \left(-A^T\right)\right)$. Note also from Proposition 2.1, that $X \in \overline{C}_A$ iff $v(X) \in \text{Im}\left(A \oplus \left(-A^T\right)\right)$. Thus, $X \in \overline{C}_A$ iff $v(X)$ is a solution to the algebraic compatible linear system:

$$\left(A \oplus \left(-A^T\right)\right) v(X) = v(M)$$

for any $M(\neq 0) \in \mathbf{R}^{n \times n}$ such that:

$$\text{rank}\left(A \oplus \left(-A^T\right)\right) = \text{rank}\left(A \oplus \left(-A^T\right), v(M)\right) = n^2 - \overline{v}(0)$$

From Theorem 2.3, the nullity and the rank of $A \oplus \left(-A^T\right)$ are, respectively, $\dim \text{Ker}\left(A \oplus \left(-A^T\right)\right) = \overline{v}(0)$ $\text{rank}\left(A \oplus \left(-A^T\right)\right) = n^2 - \overline{v}(0)$. Therefore, there exist permutation matrices $E, F \in \mathbf{R}^{n^2 \times n^2}$ such that there exists an equivalence transformation:

$$A \oplus \left(-A^T\right) \approx \overline{A} := E\left(A \oplus \left(-A^T\right)\right) F = \text{Block matrix}\left(\overline{A}_{ij}; i, j \in \overline{2}\right)$$

such that $\overline{A}_{11}$ is square non-singular and of order $\overline{v}_0$. Define $\overline{M} = EMF \approx M(\neq 0) \in \mathbf{R}^{n \times n}$. Then, the linear algebraic systems $\left(A \oplus \left(-A^T\right)\right) v(X) = v(M)$, and

$$E\left(A \oplus \left(-A^T\right)\right) F v(\overline{X}) = \begin{bmatrix} \overline{A}_{11} & \overline{A}_{12} \\ \overline{A}_{21} & \overline{A}_{22} \end{bmatrix} \begin{bmatrix} v(\overline{X}_1) \\ v(\overline{X}_2) \end{bmatrix} = \begin{bmatrix} v(\overline{M}_1) \\ v(\overline{M}_2) \end{bmatrix}$$

$$\Leftrightarrow \begin{matrix} v(\overline{X}_1) = \overline{A}_{11}^{-1}\left(v(\overline{M}_1) - \overline{A}_{12} v(\overline{X}_2)\right) \\ \left(\overline{A}_{22} - \overline{A}_{21} \overline{A}_{11}^{-1} \overline{A}_{12}\right) v(\overline{X}_2) = v(\overline{M}_2) - \overline{A}_{21} \overline{A}_{11}^{-1} v(\overline{M}_1) \end{matrix} \qquad (2.9)$$

are identical if $\overline{X}$ and $\overline{M}$ are defined according to $v(X) = Fv(\overline{X})$ and $v(\overline{M}) = Ev(M)$. As a result, Properties (i)-(ii) follow directly from (2.9) for $M = \overline{M} = 0$ and for any $M$ satisfying $\text{rank}\left(A \oplus \left(-A^T\right)\right) = \text{rank}\left(A \oplus \left(-A^T\right), v(M)\right) = n^2 - \overline{v}(0)$, respectively. □

## 3. Results concerning sets of pair-wise commuting matrices

Consider the following sets:

(1) A set of nonzero $p \geq 2$ distinct pair-wise commuting matrices $\mathbf{A}_C := \left\{A_i \in \mathbf{R}^{n \times n} : [A_i, A_j] = 0; \forall i, j \in \overline{p}\right\}$

(2) The set of matrices $MC_{\mathbf{A}_C} := \left\{X \in \mathbf{R}^{n \times n} : [X, A_i] = 0; \forall A_i \in \mathbf{A}_C\right\}$ which commute with the set $\mathbf{A}_C$ of pair-wise commuting matrices.

(3) A set of matrices $C_{\mathbf{A}} := \left\{X \in \mathbf{R}^{n \times n} : [X, A_i] = 0; \forall A_i \in \mathbf{A}\right\}$ which commute with a given set of nonzero p matrices $\mathbf{A} := \left\{A_i \in \mathbf{R}^{n \times n}; \forall i \in \overline{p}\right\}$ which are not necessarily pair-wise commuting. The complementary sets of $MC_{\mathbf{A}_C}$ and $C_{\mathbf{A}}$ are $\overline{MC_{\mathbf{A}_C}}$ and $\overline{C_{\mathbf{A}}}$, respectively, so that $\mathbf{R}^{n \times n} \ni B \in \overline{MC_{\mathbf{A}_C}}$ if $B \notin MC_{\mathbf{A}_C}$ and $\mathbf{R}^{n \times n} \ni B \in \overline{C_{\mathbf{A}}}$ if $B \notin C_{\mathbf{A}}$. Note that $C_{\mathbf{A}_C} = MC_{\mathbf{A}_C}$ for a set of pair-wise commuting matrices $\mathbf{A}_C$ so that the notation $MC_{\mathbf{A}_C}$ is directly referred to a set of matrices which commute with all those in a set of pair-wise commuting matrices. The following two basic results follow concerning commutation and non-commutation of two matrices:

**Proposition 3.1.** The following properties hold:

**(i)** $A_i \in \mathbf{A}_C; \forall i \in \overline{p} \Leftrightarrow v(A_i) \in \underset{j(\neq i) \in \overline{p}}{\cap} \text{Ker}\left(A_j \oplus \left(-A_j^T\right)\right); \forall i \in \overline{p}$

$$\Leftrightarrow v(A_i) \in \underset{i+1 \leq j \leq p}{\cap} \text{Ker}\left(A_j \oplus \left(-A_j^T\right)\right); \forall i \in \overline{p}$$

**(ii)** Define

$N_i(\mathbf{A}_C) := \left[A_1^T \oplus (-A_1) \quad A_2^T \oplus (-A_2) \quad \cdots \quad A_{i-1}^T \oplus (-A_{i-1}) \quad A_{i+1}^T \oplus (-A_{i+1}) \quad \cdots \quad A_p^T \oplus (-A_p)\right]^T$

$\in \mathbf{R}^{(p-1)n^2 \times n^2}$. Then $A_i \in \mathbf{A}_C; \forall i \in \overline{p}$ iff $v(A_i) \in \text{Ker} N_i(\mathbf{A}_C); \forall i \in \overline{p}$

(iii) $MC_{\mathbf{A}_C} := \left\{ X \in \mathbf{R}^{n \times n} : v(X) \in \bigcap_{i \in \bar{p}} Ker\left(A_i \oplus \left(-A_i^T\right)\right); A_i \in \mathbf{A}_C \right\}$

$= \left\{ X \in \mathbf{R}^{n \times n} : v(X) \in Ker N(\mathbf{A}_C) \right\} \supset C_{\mathbf{A}_C} \supset \mathbf{A}_C \supset \{0\} \in \mathbf{R}^{n \times n}$

where $N(\mathbf{A}_C) := \left[ A_1^T \oplus (-A_1) \ A_2^T \oplus (-A_2) \ \cdots \ A_p^T \oplus (-A_p) \right]^T \in \mathbf{R}^{pn^2 \times n^2}, A_i \in \mathbf{A}_C$

(iv) $\overline{MC_{\mathbf{A}_C}} := \left\{ X \in \mathbf{R}^{n \times n} : v(X) \in \bigcup_{i \in \bar{p}} Im\left(A_i \oplus \left(-A_i^T\right)\right); A_i \in \mathbf{A}_C \right\}$

$= \left\{ X \in \mathbf{R}^{n \times n} : v(X) \in Im N(\mathbf{A}_C) \right\}$

(v) $C_{\mathbf{A}} := \left\{ X \in \mathbf{R}^{n \times n} : v(X) \in \bigcap_{i \in \bar{p}} Ker\left(A_i \oplus \left(-A_i^T\right)\right); A_i \in \mathbf{A} \right\}$

$= \left\{ X \in \mathbf{R}^{n \times n} : v(X) \in Ker N(\mathbf{A}) \right\}$

where $N(\mathbf{A}) := \left[ A_1^T \oplus (-A_1) \ A_2^T \oplus (-A_2) \ \cdots \ A_p^T \oplus (-A_p) \right]^T \in \mathbf{R}^{pn^2 \times n^2}, A_i \in \mathbf{A}$

(vi) $\overline{C_{\mathbf{A}}} := \left\{ X \in \mathbf{R}^{n \times n} : v(X) \in \bigcup_{i \in \bar{p}} Im\left(A_i \oplus \left(-A_i^T\right)\right); A_i \in \mathbf{A} \right\}$

$= \left\{ X \in \mathbf{R}^{n \times n} : v(X) \in Im N(\mathbf{A}) \right\}$

**Proof**: **(i)** The first part of Property (i) follows directly from Proposition 2.1 since all the matrices of $\mathbf{A}_C$ pair-wise commute and any arbitrary matrix commutes with itself ( thus j = i may be removed from the intersections of kernels of the first double sense implication). The last part of Property (i) follows from the anti-symmetric property of the commutator $[A_i, A_j] = [A_j, A_i] = 0; \forall A_i, A_j \in \mathbf{A}_C$ what implies $A_i \in \mathbf{A}_C; \forall i \in \bar{p} \Leftrightarrow v(A_i) \in \bigcap_{i+1 \leq j \leq p} Ker\left(A_j \oplus \left(-A_j^T\right)\right); \forall A_i, A_j \in \mathbf{A}_C$.

**(ii)** It follows from its equivalence with Property (i) since $Ker N_i(\mathbf{A}_C) \equiv \bigcap_{j(\neq i) \in \bar{p}} Ker\left(A_j \oplus \left(-A_j^T\right)\right)$.

**(iii)** Property (iii) is similar to Property (i) for the whose set $M_{\mathbf{A}_C}$ of matrices which commute with the set $\mathbf{A}_C$ so that it contains $\mathbf{A}_C$ and, furthermore, $Ker N(\mathbf{A}_C) \equiv \bigcap_{i \in \bar{p}} Ker\left(A_i \oplus \left(-A_i^T\right)\right)$.

**(iv)** It follows from $\bigcup_{j\in\overline{p}} \text{Im}\left(A_j \oplus \left(-A_j^T\right)\right) = \overline{\bigcap_{j\in\overline{p}} \text{Ker}\left(A_j \oplus \left(-A_j^T\right)\right)}$ ; $A_j \in \mathbf{A}_C$ and $\mathbf{R}^{n^2} \ni 0 \in \text{Ker}\left(A_j \oplus \left(-A_j^T\right)\right) \cap \text{Im}\left(A_j \oplus \left(-A_j^T\right)\right)$ but $\mathbf{R}^{n\times n} \ni X = 0$ commutes with any matrix in $\mathbf{R}^{n\times n}$ so that $\mathbf{R}^{n\times n} \ni 0 \in MC_{\mathbf{A}_C} \Leftrightarrow \mathbf{R}^{n\times n} \ni 0 \notin \overline{MC_{\mathbf{A}_C}}$ for any given $\mathbf{A}_C$.

**(v)-(vi)** are similar to (ii)-(iv) except that the members of $\mathbf{A}$ do not necessarily commute. □

Concerning Proposition 3.1 (v)-(vi), note that if $X \in \overline{C}_{\mathbf{A}}$ then $X \neq 0$ since $\mathbf{R}^{n\times n} \ni 0 \in C_{\mathbf{A}}$. The following result is related to the rank defectiveness of the matrix $N(\mathbf{A}_C)$ and any of their sub-matrices since $\mathbf{A}_C$ is a set of pair-wise commuting matrices:

**Proposition 3.2**. The following properties hold:

$$n^2 > \text{rank } N(\mathbf{A}_C) \geq \text{rank } N_i(\mathbf{A}_C) \geq \text{rank}\left(A_j \oplus \left(-A_j^T\right)\right); \forall A_j \in \mathbf{A}_C; \forall i,j \in \overline{p}$$

and, equivalently,

$$\det\left(N^T(\mathbf{A}_C)N(\mathbf{A}_C)\right) = \det\left(N_i^T(\mathbf{A}_C)N_i(\mathbf{A}_C)\right) = \det\left(A_j \oplus \left(-A_j^T\right)\right) = 0;$$

$\forall A_j \in \mathbf{A}_C; \forall i,j \in \overline{n}$.

**Proof**: It is a direct consequence from Proposition 3.1 (i) –(ii) since the existence of nonzero pair-wise commuting matrices (all the members of $\mathbf{A}_C$) implies that the above matrices $N(\mathbf{A}_C), N_i(\mathbf{A}_C), A_j \oplus \left(-A_j^T\right)$ are all rank defective and have at least identical number of rows than that of columns. Therefore, the square matrices $N^T(\mathbf{A}_C)N(\mathbf{A}_C), N_i^T(\mathbf{A}_C)N_i(\mathbf{A}_C)$ and $A_j \oplus \left(-A_j^T\right)$ are all singular. □

Results related to sufficient conditions for a set of matrices to pair-wise commute are abundant in the literature. For instance, diagonal matrices are pair-wise commuting. Any sets of matrices taking via multiplication by real scalars with any arbitrary matrix consist of pair-wise commuting matrices. Any set of matrices obtained by linear combinations of one of the above sets consist also of pair-wise commuting matrices. Any matrix commutes with any of its matrix functions etc. In the following, we discuss a simple, although restrictive, sufficient condition for rank defectiveness of $N(\mathbf{A})$ of some set $\mathbf{A}$ of p square real n- matrices which may be useful as a test to elucidate the existence of a nonzero n- square matrix which commutes with all matrices in this set. Another useful test obtained from the following result relies on a necessary condition to elucidate if the given set consists of pair-wise commuting matrices.

**Theorem 3.3:** Consider any arbitrary set of nonzero n-square real matrices $\mathbf{A} := \{A_1, A_2, ..., A_p\}$ for any integer $p \geq 1$ and define matrices:

$$N_i(\mathbf{A}) := \begin{bmatrix} A_1^T \oplus (-A_1) & A_2^T \oplus (-A_2) & \cdots & A_{i-1}^T \oplus (-A_{i-1}) & A_{i+1}^T \oplus (-A_{i+1}) & \cdots & A_p^T \oplus (-A_p) \end{bmatrix}^T$$

$$N(\mathbf{A}) := \begin{bmatrix} A_1^T \oplus (-A_1) & A_2^T \oplus (-A_2) & \cdots & A_p^T \oplus (-A_p) \end{bmatrix}^T$$

Then, the following properties hold:

**(i)** $\mathrm{rank}(A_i \oplus (-A_i)) \le \mathrm{rank}\, N_i(\mathbf{A}) \le \mathrm{rank}\, N(\mathbf{A}) < n^2;\ \forall i \in \bar{p}$.

**(ii)** $\bigcap_{i \in \bar{p}} \mathrm{Ker}(A_i \oplus (-A_i^T)) \ne \{0\}$ so that:

$$\exists X(\ne 0) \in C_\mathbf{A},\ X \in C_\mathbf{A} \Leftrightarrow v(X) \in \bigcap_{i \in \bar{p}} \mathrm{Ker}(A_i \oplus (-A_i^T))\ \text{and}\ X \in \bar{C}_\mathbf{A} \Leftrightarrow v(X) \in \bigcup_{i \in \bar{p}} \mathrm{Im}(A_i \oplus (-A_i^T))$$

**(iii)** If $\mathbf{A} = \mathbf{A}_C$ is a set of pair-wise commuting matrices then

$$v(A_i) \in \bigcap_{j \in \bar{p} \setminus \bar{i}} \mathrm{Ker}(A_j \oplus (-A_j^T))\ ;\forall i \in \bar{p} \Leftrightarrow v(A_i) \in \bigcap_{i \in \bar{p}} \mathrm{Ker}(A_i \oplus (-A_i^T))\ ;\forall i \in \bar{p}$$

$$\Leftrightarrow v(A_i) \in \bigcap_{i \in \bar{p} \setminus \{i\}} \mathrm{Ker}(A_i \oplus (-A_i^T))\ ;\forall i \in \bar{p}$$

**(iv)** $M\mathbf{A}_C := \left\{ X \in \mathbf{R}^{n \times n} : v(X) \bigcap_{i \in \bar{p}} \mathrm{Ker}(A_i \oplus (-A_i^T)), \forall A_i \in \mathbf{A}_C \right\} \supset \mathbf{A}_C \cup \{0\} \in \mathbf{R}^{n \times n}$

with the above set inclusión being proper

**Proof**: **(i)** Any nonzero matrix $\Lambda = \mathrm{diag}(\lambda\ \lambda\ \ldots\ \lambda)$, $\lambda(\ne 0) \in \mathbf{R}$ is such that $\Lambda(\ne 0) \in C_{A_i}$ $(\forall i \in \bar{p})$ so that $\Lambda \in C_\mathbf{A}$. Thus, $0 \ne v(\Lambda) \in \mathrm{Ker}\, N(\mathbf{A}) \Leftrightarrow n^2 > \mathrm{rank}\, N(\mathbf{A}) \ge \mathrm{rank}\, N_i(\mathbf{A}) \ge \mathrm{rank}(A_i \oplus (-A_i))$; $\forall i \in \bar{p}$ and any given set $\mathbf{A}$. Property (i) has been proved.

**(ii)** The first part follows by contradiction. Assume $\bigcap_{i \in \bar{p}} \mathrm{Ker}(A_i \oplus (-A_i^T)) = \{0\}$ then $0 \ne v(\Lambda) \notin \mathrm{Ker}\, N(\mathbf{A})$ so that $\Lambda = \mathrm{diag}(\lambda\ \lambda\ \ldots\ \lambda) \notin C_\mathbf{A}$, for any $\lambda(\ne 0) \in \mathbf{R}$ what contradicts (i). Also, $X \in C_{A_i} \Leftrightarrow v(X) \in \mathrm{Ker}(A_i \oplus (-A_i^T))$; $\forall i \in \bar{p}$ so that $X \in C_\mathbf{A} \Leftrightarrow v(X) \in \bigcap_{i \in \bar{p}} \mathrm{Ker}(A_i \oplus (-A_i^T))$ what is equivalent to its contrapositive logic proposition $X \in \bar{C}_\mathbf{A} \Leftrightarrow v(X) \in \bigcup_{i \in \bar{p}} \mathrm{Im}(A_i \oplus (-A_i^T))$.

**(iii)** $\mathbf{A} = \mathbf{A}_C \Leftrightarrow A_i \in C_{A_j}; \forall j(\ne i) \in \bar{p}, \forall i \in \bar{p}$

$\Leftrightarrow A_i \in C_{A_j}; \forall j, i \in \bar{p}$ since $A_i \in C_{A_i}$ ; $\forall i \in \bar{p}$

$v(A_i) \in \bigcap_{i \in \bar{p}} \mathrm{Ker}(A_j \oplus (-A_j^T))\ ;\forall i \in \bar{p} \Leftrightarrow v(A_i) \in \bigcap_{i \in \bar{p} \setminus \{i\}} \mathrm{Ker}(A_j \oplus (-A_j^T))\ ;\forall i \in \bar{p}$

On the other hand,

$(v(A_i) \in \bigcap_{j \in \bar{p} \setminus i} \mathrm{Ker}(A_j \oplus (-A_j^T)) \Leftrightarrow v(A_i) \in C_{A_j}; \forall j \in \bar{p})$ for any $i\,(<p) \in \bar{p}$.

This assumption implies directly that:

$$v(A_i) \in C_{A_j}; \forall j \in \bar{p} \,\wedge\, v(A_{i+1}) \in \bigcap_{j \in \overline{i+1}} C_{A_j} \text{ for any } i\,(<p) \in \bar{p}$$

which together with $v(A_{i+1}) \in \bigcap_{j \in \bar{p} \setminus \overline{i+1}} \mathrm{Ker}(A_j \oplus (-A_j^T))$ implies that $v(A_{i+1}) \in C_{A_j}; \forall j \in \bar{p}$

$$\Rightarrow (v(A_{i+1}) \in \bigcap_{j \in \bar{p} \setminus \overline{i+1}} \mathrm{Ker}(A_j \oplus (-A_j^T)) \Leftrightarrow \,) \text{ for } (i+1) \in \bar{p}$$

Thus, it follows by complete induction that $A = A_C \Leftrightarrow v(A_i) \in \bigcap_{i \in \bar{p} \setminus \{i\}} \mathrm{Ker}(A_j \oplus (-A_j^T))$; $\forall i \in \bar{p}$ and Property (iii) has been proved.

**(iv)** The definition of $MA_C$ follows from Property (iii) in order to guarantee that $[X, A_i] = 0$; $\forall A_i \in A$. The fact that such a set contains properly $A_C \cup \{0\}$ follows directly from $\mathbf{R}^{n \times n} \ni \Lambda = \mathrm{diag}(\lambda\,\lambda\,\cdots\,\lambda)(\in MC_{A_C}) \neq A_C \cup \{0\}$ for any $\mathbf{R} \ni \lambda \neq 0$. □

Note that Theorem 3.3 (ii) extends Proposition 3.1 (v) since it is proved that $C_A \setminus \{0\} \neq \varnothing$ because all nonzero $\mathbf{R}^{n \times n} \ni \Lambda = \mathrm{diag}(\lambda\,\lambda\,\cdots\,\lambda) \in C_A$ for any $\mathbf{R} \ni \lambda \neq 0$ and any set of matrices $A$. Note that Theorem 3.3 (iii) establishes that $v(A_i) \in \bigcap_{i \in \bar{p} \setminus \{i\}} \mathrm{Ker}(A_j \oplus (-A_j^T))$; $\forall i \in \bar{p}$ is a necessary and sufficient condition for the set to be a set of commuting matrices $A$ being simpler to test (by taking advantage of the symmetry property of the commutators) than the equivalent condition $v(A_i) \in \bigcap_{i \in \bar{p}} \mathrm{Ker}(A_j \oplus (-A_j^T))$; $\forall i \in \bar{p}$. Further results about pair-wise commuting matrices or the existence of nonzero commuting matrices with a given set are obtained in the subsequent result based on the Kronecker sum of relevant Jordan canonical forms:

**Theorem 3.4.** The following properties hold for any given set of n-square real matrices $A = \{A_1, A_2, \ldots, A_p\}$:

**(i)** The set $C_A$ of matrices $X \in \mathbf{R}^{n \times n}$ which commute with all matrices in $A$ is defined by:

$$C_A := \left\{ X \in \mathbf{R}^{n \times n} : v(X) \in \bigcap_{i=1}^{p} \left( \mathrm{Ker}\left[ (J_{A_i} \oplus (-J_{A_i}^T))(P_i^{-1} \otimes P_i^{-T}) \right] \right) \right\} \quad (3.1)$$

$$= \left\{ X \in \mathbf{R}^{n \times n} : v(X) \in \bigcap_{i=1}^{p} \left( \mathrm{Im}\left( (P_i \otimes P_i^{-1})(Y_i) \right) \right) \wedge Y_i \in \mathrm{Ker}(J_{A_i} \oplus (-J_{A_i}^T)); \forall i \in \bar{p} \right\}$$
(3.2)

$$= \left\{ X \in \mathbf{R}^{n \times n} : v(X) \in \bigcap_{i=1}^{p} \left( \mathrm{Im}\left( (P_i \otimes P_i^{-1})(Y) \right) \right), Y \in \bigcap_{i=1}^{p} \left( \mathrm{Ker}(J_{A_i} \oplus (-J_{A_i}^T)) \right) \right\}$$
(3.3)

where $P_i \in \mathbf{R}^{n \times n}$ is a non-singular transformation matrix such that $A_i \sim J_{A_i} = P_i^{-1} A_i P_i$, $J_{A_i}$ being the Jordan canonical form of $A_i$.

**(ii)** $\dim \mathrm{span}\{v(X): X \in C_\mathbf{A}\} \leq \min_{i \in \bar{p}} \dim\left(\mathrm{Ker}\left(J_{A_i} \oplus \left(-J_{A_i}^T\right)\right)\right)$

$$= \min_{i \in \bar{p}} \bar{v}_i(0) = \min_{i \in \bar{p}} \left(\sum_{j=1}^{\rho_i} v_{ij}^2\right) \leq \min_{i \in \bar{p}} \left(\sum_{i=1}^{\rho_i} \mu_{ij}^2\right) \leq \min_{i \in \bar{p}} \left(\bar{\mu}_i(0)\right)$$

where $\bar{v}_i(0)$ and $v_{ij}$ are, respectively, the geometric multiplicities of $0 \in \sigma\left(A_i \oplus \left(-A_i^T\right)\right)$ and $\lambda_{ij} \in \sigma(A_i)$ and $\bar{\mu}_i(0)$ and $\mu_{ij}$ are, respectively, the algebraic multiplicities of $0 \in \sigma\left(A_i \oplus \left(-A_i^T\right)\right)$ and $\lambda_{ij} \in \sigma(A_i)$; $\forall j \in \rho_i$ (the number of distinct eigenvalues of $A_i$), $\forall i \in \bar{p}$.

**(iii)** The set $\mathbf{A}$ consists of pair-wise commuting matrices, namely $C_\mathbf{A} = MC_\mathbf{A}$, iff

$$v(A_j) \in \bigcap_{i(\neq j)=1}^{p} \left(\mathrm{Ker}\left[\left(J_{A_i} \oplus \left(-J_{A_i}^T\right)\right)\left(P_i^{-1} \otimes P_i^{-T}\right)\right]\right); \quad \forall j \in \bar{p}.$$

Equivalent conditions follow from the second and third equivalent definitions of $C_\mathbf{A}$ in Property (i).

**Proof**: If $A_i \sim J_{A_i} = P_i^{-1} A_i P_i$, with $J_{A_i}$ being the Jordan canonical form of $A_i$ then $A_i \oplus \left(-A_i^T\right) \sim J_{A_i} \oplus \left(-J_{A_i}^T\right) = T_i^{-1}\left(A_i \oplus \left(-A_i^T\right)\right) T_i$ with $T_i = P_i \otimes P_i^T \in \mathbf{R}^{n^2} \times \mathbf{R}^{n^2}$ (see proof of Theorem 2.3) being non-singular; $\forall i \in \bar{p}$. Thus, $\left(A_i \oplus \left(-A_i^T\right)\right) = T_i \left(J_{A_i} \oplus \left(-J_{A_i}^T\right)\right) T_i^{-1}$ so that:

$$N(\mathbf{A}) = \left[A_1^T \oplus (-A_1) \; A_2^T \oplus (-A_2) \; \cdots \; A_p^T \oplus (-A_p)\right]^T = \left[I_{pn^2} \; -U\right]\left[W_1^T \; W_2^T\right]^T = \mathbf{TJT_a}$$

(3.4)

where

$$\mathbf{T} := \mathrm{Block\,Diag}\left[T_1 \; T_2 \ldots T_p\right] \in \mathbf{R}^{pn^2 \times pn^2}; \quad \mathbf{T_a} := \left[T_1^{-T} \; T_2^{-T} \; \cdots \; T_p^{-T}\right]^T \in \mathbf{R}^{pn^2 \times n^2}$$

(3.5)

$$\mathbf{J} := \mathrm{Block\,Diag}\left[J_{A_1} \oplus \left(-J_{A_1}^T\right) \; J_{A_2} \oplus \left(-J_{A_2}^T\right) \cdots J_{A_p} \oplus \left(-J_{A_p}^T\right)\right] \in \mathbf{R}^{pn^2 \times pn^2}$$

(3.6)

Then,

$$\mathrm{Ker}\,N(\mathbf{A}) = \bigcap_{i=1}^{p} \mathrm{Ker}\left[\left(A_i \oplus \left(-A_i^T\right)\right)\right] = \mathrm{Ker}(\mathbf{JT_a}) \equiv \bigcap_{i=1}^{p} \left(\mathrm{Ker}\left[\left(J_{A_i} \oplus \left(-J_{A_i}^T\right)\right)\left(P_i^{-1} \otimes P_i^{-T}\right)\right]\right)$$

since $\mathbf{T}$ is non-singular. Thus, $\forall X \in \mathrm{Dom}(\mathbf{A}) \subset \mathbf{R}^{n^2}$:

$X \in C_\mathbf{A} \Leftrightarrow v(X) \in \mathrm{Ker}\,N(\mathbf{A})$

$$\Leftrightarrow v(X) \in \left( \bigcap_{i=1}^{p} \operatorname{Ker}\left[ \left(A_i \oplus \left(-A_i^T\right)\right) \right] \right) \Leftrightarrow v(X) \in \left( \bigcap_{i=1}^{p} \operatorname{Ker}\left[ \left(J_{A_i} \oplus \left(-J_{A_i}^T\right)\right)\left(P_i^{-1} \otimes P_i^{-T}\right) \right] \right)$$

$$\Leftrightarrow v(X) \in \operatorname{Im}\left( \left(P_i \otimes P_i^{-1}\right) \left( \operatorname{Ker}\left(J_{A_i} \oplus \left(-J_{A_i}^T\right)\right) \right) \right); \forall i \in \overline{p}$$

$$\Leftrightarrow v(X) \in \bigcap_{i=1}^{p} \left( \operatorname{Im}\left( \left(P_i \otimes P_i^{-1}\right)(Y) \right) \right) \Leftrightarrow v(X) \in \bigcap_{i=1}^{p} \left( \operatorname{Im}\left( \left(P_i \otimes P_i^{-1}\right)(Y_i) \right) \right)$$

where $Y_i \in \operatorname{Ker}\left(J_{A_i} \oplus \left(-J_{A_i}^T\right)\right)$; $\forall i \in \overline{p}$ and $Y \in \left( \bigcap_{i=1}^{p} \left( \operatorname{Ker}\left(J_{A_i} \oplus \left(-J_{A_i}^T\right)\right) \right) \right)$. Property (i) has been proved. The first inequality of Property (ii) follows directly from Property (i). The results of equalities and inequalities in the second line of Property (ii) follow by the first inequality by taking into account and Theorem 2.3. Property (iii) follows from the proved equivalent definitions of $C_A$ in Property (i) by taking into account that $[A_j, A_j] = 0$; $\forall j \in \overline{p}$ so that:

$$v(A_j) \in \bigcap_{i=1}^{p} \left( \operatorname{Ker}\left[ \left(J_{A_i} \oplus \left(-J_{A_i}^T\right)\right)\left(P_i^{-1} \otimes P_i^{-T}\right) \right] \right)$$

$$\Leftrightarrow v(A_j) \in \bigcap_{i(\neq j)=1}^{p} \left( \operatorname{Ker}\left[ \left(J_{A_i} \oplus \left(-J_{A_i}^T\right)\right)\left(P_i^{-1} \otimes P_i^{-T}\right) \right] \right); \forall j \in \overline{p} \qquad \square$$

Theorem 3.3 are concerned with $MC_A \neq \{0\} \in \mathbf{R}^{n \times n}$ for an arbitrary set of real square matrices A and for a pair-wise-commuting set, respectively.

## 4. Further results and extensions

The extensions of the results for commutation of complex matrices is direct in several ways. It is first possible to decompose the commutator in its real and imaginary part and then apply the results of Sections 2-3 for real matrices to both parts as follows. Let $A = A_{re} + \mathbf{i} A_{im}$ and $B = B_{re} + \mathbf{i} B_{im}$ be complex matrices in $\mathbf{C}^{n \times n}$ with $A_{re}$, and $B_{re}$ being their respective real parts, and $A_{im}$ and $B_{im}$, all in $\mathbf{R}^{n \times n}$ their respective imaginary parts and $\mathbf{i} = \sqrt{-1}$ is the imaginary complex unity. Direct computations with the commutator of A and B yield:

$$[A, B] = \left([A_{re}, B_{re}] - [A_{im}, B_{im}]\right) + \mathbf{i}\left([A_{im}, B_{re}] + [A_{re}, B_{im}]\right) \qquad (4.1)$$

The following three results are direct and allow to reduce the problem of commutation of a pair of complex matrices to the discussion of four real commutators:

**Proposition 4.1.** $B \in C_A \Leftrightarrow \left( \left([A_{re}, B_{re}] = [A_{im}, B_{im}]\right) \wedge \left([A_{im}, B_{re}] = [B_{im}, A_{re}]\right) \right)$.

**Proposition 4.2.** $\left(B_{re} \in \left(C_{A_{re}} \cap C_{A_{im}}\right) \wedge B_{im} \in \left(C_{A_{im}} \cap C_{A_{re}}\right)\right) \Rightarrow B \in C_A$.

**Proposition 4.3.** $\left(A_{re} \in \left(C_{B_{re}} \cap C_{B_{im}}\right) \wedge A_{im} \in \left(C_{B_{im}} \cap C_{B_{re}}\right)\right) \Rightarrow B \in C_A$.

**Proofs**: Proposition 4.1 follow by inspection of (4.1). Proposition 4.2 implies that Proposition 4.1 holds with the four involved commutators being zero. Then the left condition of Proposition 4.2 implies that $B \in C_A$, from Proposition 4.1, so that Proposition 4.2 holds. Proposition 4.3 is equivalent to Proposition 4.2. □

Proposition 4.1 yields to the subsequent result

**Theorem 4.4.** The following properties hold:

**(i)** Assume that the matrices A and $B_{re}$ are given. Then, $B \in C_A$ iff $B_{im}$ satisfies the linear algebraic equation:

$$\begin{bmatrix} A_{re} \oplus \left(-A_{re}^T\right) \\ A_{im} \oplus \left(-A_{im}^T\right) \end{bmatrix} v(B_{re}) = \begin{bmatrix} A_{im} \oplus \left(-A_{im}^T\right) \\ A_{re} \oplus \left(-A_{re}^T\right) \end{bmatrix} v(B_{im}) \quad (4.2)$$

for which a necessary condition is:

$$\operatorname{rank}\begin{bmatrix} A_{im} \oplus \left(-A_{im}^T\right) \\ A_{re} \oplus \left(-A_{re}^T\right) \end{bmatrix} = \operatorname{rank}\begin{bmatrix} A_{im} \oplus \left(-A_{im}^T\right) & \left( A_{re} \oplus \left(-A_{re}^T\right) \right) v(B_{re}) \\ A_{re} \oplus \left(-A_{re}^T\right) & \left( A_{im} \oplus \left(-A_{im}^T\right) \right) \end{bmatrix} \quad (4.3)$$

**(ii)** Assume that the matrices A and $B_{ime}$ are given. Then, $B \in C_A$ iff $B_{re}$ satisfies (4.2) for which a necessary condition is:

$$\operatorname{rank}\begin{bmatrix} A_{re} \oplus \left(-A_{re}^T\right) \\ A_{im} \oplus \left(-A_{im}^T\right) \end{bmatrix} = \operatorname{rank}\begin{bmatrix} A_{re} \oplus \left(-A_{re}^T\right) & \left( A_{im} \oplus \left(-A_{im}^T\right) \right) v(B_{im}) \\ A_{im} \oplus \left(-A_{im}^T\right) & \left( A_{re} \oplus \left(-A_{re}^T\right) \right) \end{bmatrix}$$

**(iii)** Also, $\exists B \neq 0$ such that $B \in C_A$ with $B_{re} = 0$ and $\exists B \neq 0$ such that $B \in C_A$ with $B_{im} = 0$

Also, $\exists B_{im} \neq 0$ Instead, if

**Proof**: **(i)** Eqn. 4.2 is a re-arrangement in an equivalent algebraic system of Proposition 4.1 in the unknown $v(B_{im})$ for given A and $B_{re}$. The system is compatible if (4.2) holds from the Kronecker-Capelli theorem. The proof of Property (ii) is similar to that of (i) with the appropriate interchange of roles of $B_{re}$ and $B_{im}$.

**(iii)** Since $\operatorname{rank}\begin{bmatrix} A_{im} \oplus \left(-A_{im}^T\right) \\ A_{re} \oplus \left(-A_{re}^T\right) \end{bmatrix} < n^2$ from Theorem 3.3 (i) then $0 \neq B = B_{re} \in C_A$ iff $B_{re} \in C_{A_{re}} \cap C_{A_{im}} (\neq \emptyset) \subset C_A$. The same proof follows for $0 \neq B = B_{im} \in C_A$ since

$$\operatorname{rank}\begin{bmatrix} A_{re} \oplus \left(-A_{re}^T\right) \\ A_{im} \oplus \left(-A_{im}^T\right) \end{bmatrix} = \operatorname{rank}\begin{bmatrix} A_{im} \oplus \left(-A_{im}^T\right) \\ A_{re} \oplus \left(-A_{re}^T\right) \end{bmatrix} < n^2. \quad □$$

A more general result than Theorem 4.4 is the following:

**Theorem 4.5**. The following properties hold:

**(i)** $B \in C_A \cap \mathbf{C}^{n \times n}$ iff $v(B)$ is a solution to the linear algebraic system :

$$\begin{bmatrix} A_{re} \oplus (-A_{re}^T) & (-A_{im}) \oplus (A_{im}^T) \\ A_{im} \oplus (-A_{im}^T) & (-A_{re}) \oplus (A_{re}^T) \end{bmatrix} \begin{bmatrix} v(B_{re}) \\ v(B_{im}) \end{bmatrix} = 0 \qquad (4.4)$$

Nonzero solutions $B \in C_A$, satisfying $\begin{bmatrix} v(B_{re}) \\ v(B_{im}) \end{bmatrix} \in \mathrm{Ker} \begin{bmatrix} A_{re} \oplus (-A_{re}^T) & (-A_{im}) \oplus (A_{im}^T) \\ A_{im} \oplus (-A_{im}^T) & (-A_{re}) \oplus (A_{re}^T) \end{bmatrix}$,

always exist since

$$\mathrm{Ker} \begin{bmatrix} A_{re} \oplus (-A_{re}^T) & (-A_{im}) \oplus (A_{im}^T) \\ A_{im} \oplus (-A_{im}^T) & (-A_{re}) \oplus (A_{re}^T) \end{bmatrix} \neq \{0\} \in \mathbf{R}^{2n^2} \text{, and equivalently, since}$$

$$\mathrm{rank} \begin{bmatrix} A_{re} \oplus (-A_{re}^T) & (-A_{im}) \oplus (A_{im}^T) \\ A_{im} \oplus (-A_{im}^T) & (-A_{re}) \oplus (A_{re}^T) \end{bmatrix} < 2n^2 \qquad (4.5)$$

**(ii)** Property (ii) is equivalent to

$$B \in C_A \Leftrightarrow \left( A \oplus (-A^*) \right) v(B) = 0 \qquad (4.6)$$

which has always nonzero solutions since $\left( A \oplus (-A^*) \right) < n^2$

**Proof**: **(i)** It follows in the same way as that of Theorem 4.4 by rewriting the algebraic system (4.3) in the form (4.4) which has nonzero solutions if (4.5) holds. But (4.5) always holds since $B = A \in C_A \cap \mathbf{C}^{n \times n}$ is nonzero if $A$ is nonzero and if $A = 0 \in \mathbf{C}^{n \times n}$ then $C_A = \mathbf{C}^{n \times n}$.

**(ii)** Direct calculations yield the equivalence of (4.4) with the separation into real and imaginary parts of the subsequent algebraic system:

$$\left( A \otimes I_n - I_n \otimes A^* \right) v(B) = \left[ \left( A_{re} + \mathbf{i} A_{im} \right) \otimes I_n - I_n \otimes \left( A_{re}^T - \mathbf{i} A_{im}^T \right) \right] \left( v(B_{re}) + \mathbf{i} v(B_{im}) \right) = 0$$

which is always solvable with a nonzero solution (i.e. compatible) since $\mathrm{rank}\left( A \otimes I_n - I_n \otimes A^* \right) < n^2$ ( otherwise , $A(\neq 0) \in C_A$ ). □

The various results of Section 3 for a set of distinct complex matrices to pair-wise commute and for characterizing the set of complex matrices which commute with those in a given set may be discussed by more general algebraic systems like the above one with four block matrices

$$\begin{bmatrix} A_{jre} \oplus \left(-A_{2re}^T\right) & \left(-A_{jim}\right) \oplus \left(A_{jim}^T\right) \\ A_{jim} \oplus \left(-A_{2im}^T\right) & \left(-A_{j2re}\right) \oplus \left(A_{jre}^T\right) \end{bmatrix}$$ for each $j \in \bar{p}$ in the whole algebraic system. Theorem 4.5 extends directly for sets of complex matrices commuting with a given one and complex matrices commuting with a set of commuting complex matrices as follows:

**Theorem 4.6**. The following properties hold:

**(i)** Consider the sets of nonzero distinct complex matrices $\mathbf{A} := \{A_i \in \mathbf{C}^{n \times n} : i \in \bar{p}\}$ and $C_\mathbf{A} := \{X \in \mathbf{C}^{n \times n} : [X, A_i] = 0; A_i \in \mathbf{A}, \forall i \in \bar{p}\}$ for $p \geq 2$. Thus, $C_\mathbf{A} \ni X = X_{re} + iX_{re}$ iff

$$\begin{bmatrix} A_{1re} \oplus \left(-A_{1re}^T\right) & \left(-A_{1im}\right) \oplus \left(A_{1im}^T\right) \\ A_{1im} \oplus \left(-A_{1im}^T\right) & \left(-A_{1re}\right) \oplus \left(A_{1re}^T\right) \\ A_{2re} \oplus \left(-A_{2re}^T\right) & \left(-A_{2im}\right) \oplus \left(A_{2im}^T\right) \\ A_{2im} \oplus \left(-A_{2im}^T\right) & \left(-A_{2re}\right) \oplus \left(A_{2re}^T\right) \\ \vdots & \\ A_{pre} \oplus \left(-A_{pre}^T\right) & \left(-A_{pim}\right) \oplus \left(A_{pim}^T\right) \\ A_{pim} \oplus \left(-A_{pim}^T\right) & \left(-A_{pre}\right) \oplus \left(A_{pre}^T\right) \end{bmatrix} \begin{bmatrix} v(X_{re}) \\ v(X_{im}) \end{bmatrix} = 0 \qquad (4.7)$$

and a nonzero solution $X \in C_\mathbf{A}$ exists since the rank of the coefficient matrix of (4.7) is less than $2n^2$.

**(ii)** Consider the sets of nonzero distinct commuting complex matrices $\mathbf{A}_C := \{A_i \in \mathbf{C}^{n \times n} : i \in \bar{p}\}$ and $MC_\mathbf{A} := \{X \in \mathbf{C}^{n \times n} : [X, A_i] = 0; A_i \in \mathbf{A}, \forall i \in \bar{p}\}$ for $p \geq 2$. Thus, $MC_\mathbf{A} \ni X = X_{re} + iX_{re}$ iff $v(X_{re})$ and $v(X_{im})$ are solutions to (4.7).

**(iii)** Properties (i) and (ii) are equivalently formulated by from the algebraic set of complex equations:

$$\left[A_1^* \oplus \left(-A_1\right) \quad A_2^* \oplus \left(-A_2\right) \quad \cdots \quad A_p^* \oplus \left(-A_p\right)\right]^* v(X) = 0 \qquad (4.8)$$

**Outline of Proof**: **(i)** It is a direct extension of Theorem 4.5 by decomposing the involved complex matrices in their real and imaginary parts since from Theorem 3.3 (i) both left block matrices in the coefficient matrix of (4.7) have rank less than $n^2$. As a result, such a coefficient matrix has rank less than $2n^2$ so that nonzero solutions exists to the algebraically compatible system of linear equations (4.7). As a result, a nonzero n-square complex commuting matrix exists.

**(ii)** It is close to that of (i) but the rank condition for compatibility of the algebraic system is not needed since the coefficient matrix of (4.7) is rank defective since $A_j \in \mathbf{A}_C \Leftrightarrow \left(v^T(A_{jre}), v^T(A_{jim})\right)^T$ is in the null space of the coefficient matrix; $\forall j \in \overline{p}$.

**(iii)** Its proof is close to that of Theorem 4.5 (ii) and it is then omitted. □

**Remark 4.7**. Note that all the proved results of Sections 2- 3 are directly extendable for complex commuting matrices, by simple replacements of transposes by conjugate transposes, without requiring a separate decomposition in real and imaginary parts as discussed in Theorem 4.5(ii) and Theorem 4.6 (iii).

□

Let $f: \mathbf{C} \to \mathbf{C}$ be an analytic function in an open set $D \supset \sigma(A)$ for some matrix $A \in \mathbf{C}^{n \times n}$ and let $p(\lambda)$ a polynomial fulfilling $p^{(i)}(\lambda_k) = f^{(i)}(\lambda_k)$; $\forall k \in \sigma(A)$, $\forall i \in \overline{m_k - 1} \cup \{0\}$; $\forall k \in \mu$ (the number of distinct elements in $\sigma(A)$, where $m_k$ is the index of $\lambda_k$, that is its multiplicity in the minimal polynomial of A. Then, f (A) is a function of a matrix A if $f(A) = p(A)$, [8]. Some results follow concerning the commutators of functions of matrices.

**Theorem 4.8**. Consider a nonzero matrix $B \in C_A \cap \mathbf{C}^{n \times n}$ for any given nonzero $A \in \mathbf{C}^{n \times n}$. Then, $f(B) \in C_A \cap \mathbf{C}^{n \times n}$, and equivalently $v(f(B)) \in \text{Ker}\left(A \oplus (-A^*)\right)$, for any function $f: \mathbf{C}^{n \times n} \to \mathbf{C}^{n \times n}$ of the matrix B.

**Proof**: For any $B \in C_A \cap \mathbf{C}^{n \times n}$:

$$[A, B] = 0 \Rightarrow (\lambda I_n - B) A = A(\lambda I_n - B); \forall \lambda \in \mathbf{C} \Rightarrow (\lambda I_n - B)^{-1} A = A(\lambda I_n - B)^{-1}; \forall \lambda \in \mathbf{C} \cap \overline{\sigma(B)}$$

(4.9)

$$\Rightarrow [A, f(B)] = A\left[\frac{1}{2\pi \mathbf{i}} \oint_C f(\lambda)(\lambda I_n - B)^{-1} d\lambda\right] = \frac{1}{2\pi \mathbf{i}} \oint_C f(\lambda)(\lambda I_n - B)^{-1} A d\lambda$$

$$= \left[\frac{1}{2\pi \mathbf{i}} \oint_C f(\lambda)(\lambda I_n - B)^{-1} d\lambda\right] A = [f(B), A] = 0$$

where C is the boundary of D and consists in a set of closed rectifiable Jordan curves which contains no point of $\sigma(A)$ since $\lambda \in \mathbf{C} \cap \overline{\sigma(A)}$ so that the identity $(\lambda I_n - B)^{-1} A = A(\lambda I_n - B)^{-1}$ is true. Then, $f(B) \in C_A \cap \mathbf{C}^{n \times n}$ has been proved. From Theorem 3.5, this is equivalent to $v(f(B)) \in \text{Ker}\left(A \oplus (-A^*)\right)$. □

The following corollaries are direct from Theorem 4.8 from the subsequent facts:

1) $A \in C_A$; $\forall A \in \mathbf{C}^{n \times n}$.
2) $[A, B] = 0 \Rightarrow [A, g(B)] = 0 \Rightarrow [f(A), g(B)]$

$$= [p(A), g(B)] = \sum_{i=0}^{\mu} \alpha_i [A^i, g(B)] = \sum_{i=0}^{\mu} \alpha_i A^{i-1}[A, g(B)] = 0 \Leftrightarrow g(B) \in C_{f(A)} \cap \mathbf{C}^{n \times n}$$

where $f(A) = p(A)$, from the definition of f being a function of the matrix A, with $p(\lambda)$ being a polynomial fulfilling $p^{(i)}(\lambda_k) = f^{(i)}(\lambda_k)$; $\forall k \in \sigma(A)$, $\forall i \in \overline{m_k - 1} \cup \{0\}$; $\forall k \in \mu$ (the number of distinct elements in $\sigma(A)$, where $m_k$ is the index of $\lambda_k$, that is its multiplicity in the minimal polynomial of A.

3) Theorem 4.8 is extendable for any countable set $\{f_i(B)\}$ of matrix functions of B.

**Corollary 4.9**. Consider a nonzero matrix $B \in C_A \cap \mathbf{C}^{n \times n}$ for any given nonzero $A \in \mathbf{C}^{n \times n}$. Then, $g(B) \in C_{f(A)} \cap \mathbf{C}^{n \times n}$ for any function $f: \mathbf{C}^{n \times n} \to \mathbf{C}^{n \times n}$ of the matrix A and any function $g: \mathbf{C}^{n \times n} \to \mathbf{C}^{n \times n}$ of the matrix B. □

**Corollary 4.10**. $f(A) \in C_A \cap \mathbf{C}^{n \times n}$, and equivalently $v(f(A)) \in \text{Ker}(A \oplus (-A^*))$, for any function $f: \mathbf{C}^{n \times n} \to \mathbf{C}^{n \times n}$ of the matrix A. □

**Corollary 4.11**. If $B \in C_A \cap \mathbf{C}^{n \times n}$ then any countable set of function matrices $\{f_i(B)\}$ is $C_A$ and in $MC_A$. □

**Corollary 4.12**. Consider any countable set of function matrices $C_F := \{f_i(A); \forall i \in \overline{p}\} \subset C_A$ for any given nonzero $A \in \mathbf{C}^{n \times n}$. Then, $\bigcap_{f_i \in C_F} \left(\text{Ker}(f_i(A) \oplus (-f_i(A^*)))\right) \supset \text{Ker}(A \oplus (-A^*))$. □

Note that matrices which commute and are simultaneously triangularizable through the same similarity transformation maintain a zero commutator after such a transformation is performed.

**Theorem 4.12**. Assume that $B \in C_A \cap \mathbf{C}^{n \times n}$, Thus, $\Lambda_B \in C_{\Lambda_A} \cap \mathbf{C}^{n \times n}$ provided that there exists a non-singular matrix $T \in \mathbf{C}^{n \times n}$ such that $\Lambda_A = T^{-1}AT$ and $\Lambda_B = T^{-1}BT$.

**Proof**: $B \in C_A \Leftrightarrow F[A, B]G = 0$; $\forall F, G \in \mathbf{C}^{n \times n}$ being non-singular. By choosing $F^{-1} = G = T$, it follows that

$$T^{-1}[A, B]T = T^{-1}A(TT^{-1})BT - T^{-1}B(TT^{-1})AT = [\Lambda_A, \Lambda_B] = 0.$$ □

A direct consequence of Theorem 4.12 is that if a set of matrices are simultaneously triangularizable to their real canonical forms by a common transformation matrix then the pair-wise commuting properties are identical to those of their respective Jordan forms.

**ACKNOWLEDGMENTS**

The author is very grateful to the Spanish Ministry of Education by its partial support of this work through project DPI2006-00714. He is also grateful to the Basque Government by its support through GIC07143-IT-269-07, SAIOTEK SPE07UN04.